\documentclass[12pt,reqno]{amsart}

\usepackage{amsthm,amsmath,amssymb,amsfonts,amscd,array,hhline,color,verbatim,latexsym,xypic}
\usepackage{fullpage} 
\newtheorem{Theorem}{Theorem}[section]
\newtheorem{Lemma}[Theorem]{Lemma}
\newtheorem{Proposition}[Theorem]{Proposition}

\newtheorem{Conjecture}[Theorem]{Conjecture} 

\newcommand{\om}{\omega}
\newcommand{\Hilb}{\mathcal H} 
\newcommand{\tGott}{\widetilde{\mathfrak G}}
\newcommand{\Gott}{\mathfrak G} 
\newcommand{\SI}{\mathfrak S} 
\newcommand{\cb}{\,) \hspace{-1.6mm} ( \,}
\newcommand{\B}{\mathbb B} 
\newcommand{\cW}{\mathcal W} 
\newcommand{\F}{\mathbb F} 
\newcommand{\tA}{\widetilde A} 
\newcommand{\tF}{\widetilde F} 
\newcommand{\He}{\text{He} }
\newcommand{\fg}{\mathfrak g} 
\newcommand{\fra}{\mathfrak a} 
\newcommand{\complex}{\mathbf C}

\renewcommand{\P}{\mathbf P}
\newcommand{\Aff}{\mathbb A}
\newcommand{\cA}{\mathcal A}
\newcommand{\cE}{\mathcal E}
\newcommand{\ua}{\underline a} 
\newcommand{\ub}{\underline b} 
\newcommand{\uu}{\mathbf u} 
\newcommand{\ux}{\mathbf x} 
\newcommand{\uy}{\mathbf y} 
\newcommand{\uz}{\mathbf z} 
\newcommand{\Proj}{\text{Proj} \,} 

\renewcommand{\le}{\leqslant}
\renewcommand{\ge}{\geqslant}

\newcommand{\lra}{\longrightarrow}
\newcommand{\la}{\leftarrow}

\newcommand{\demo}{\noindent {\sc Proof.}\;}

\def\bbone{{\mathchoice {\rm 1\mskip-4mu l} {\rm 1\mskip-4mu l}
{\rm 1\mskip-4.5mu l} {\rm 1\mskip-5mu l}}}

\begin{document}
\title[Hilbert covariants]{On Hilbert covariants} 
\author[HH]{}
\maketitle 
\centerline{Abdelmalek Abdesselam and Jaydeep Chipalkatti} 

\bigskip \bigskip 

\parbox{16cm}{ \small
{\sc Abstract:} Let $F$ denote a binary form of order $d$ over the
complex numbers. If $r$ is a divisor of $d$, then the Hilbert covariant
$\Hilb_{r,d}(F)$ vanishes exactly when $F$ is the perfect power of an
order $r$ form. In geometric terms, the coefficients of $\Hilb$ give
defining equations for the image variety $X$ of an embedding $\P^r
\hookrightarrow \P^d$. In this paper we describe a new construction of
the Hilbert covariant; and simultaneously situate it into a wider class of
covariants called the G{\"o}ttingen covariants, all of which vanish on
$X$. We prove that the ideal generated by the coefficients of $\Hilb$
defines $X$ as a scheme. Finally, we exhibit a generalisation of the
G{\"o}ttingen covariants to $n$-ary forms using the classical Clebsch transfer principle.}

\medskip \medskip 


\bigskip 

\setcounter{tocdepth}{1} 
\tableofcontents

\thispagestyle{empty} 
\section{Introduction} 
\subsection{} \label{section.intro1} 
Let 
\[ F = \sum\limits_{i=0}^d \, \binom{d}{i} \, a_i \, x_1^{d-i} \,
x_2^i,  \qquad (a_i \in \complex) \] 
denote a binary form of order\footnote{Usually $d$ would be called
  the degree of $F$, but `order' is  the common usage in classical
  invariant theory.} $d$ in the variables $\ux =
\{x_1,x_2\}$.  Its Hessian is defined to be 
\[ \He(F) = \frac{\partial^2 F}{\partial x_1^2} \, \frac{\partial^2 F}{\partial x_2^2} 
- \left(\frac{\partial^2 \, F}{\partial x_1 \, \partial x_2} \right)^2. \] 
It is well-known that  
\[ \He(F) = 0 \iff F = (p \, x_1 + q \, x_2)^d \quad \text{for some $p,q \in
  \complex$}. \] 
(The implication $\Leftarrow$ is obvious, and $\Rightarrow$ easily
follows by a simple integration -- see~\cite[Proposition~2.23]{Olver}.) 

The Hessian is a covariant of binary $d$-ics,  in the sense
that its construction commutes with a linear change of variables in
the $\ux$. More precisely, let $g = \left( \begin{array}{cc} \alpha & \gamma
  \\ \beta & \delta \end{array} \right)$ denote a complex matrix such
that $\det g = 1$. Given a binary form $A(x_1,x_2)$, write 
\[ A^g = A(\alpha \, x_1 + \beta \, x_2, \gamma \, x_1 + \delta \,
x_2). \] Then we have an identity 
\[ \He(F^g) = [ \, \He(F) \, ]^g. \] 
By definition, $\He(F)$ is of order $2d-4$ (in the $\ux$), and its coefficients are 
quadratic in the $a_i$; hence it is said to be a covariant of degree $2$ and order $2d-4$. 
\subsection{} \label{power.problem} 
Now suppose that $r$ is a divisor of $d$ (say $d = r \, \mu$), and we are looking for a
similar covariant which vanishes exactly whether $F$ is the perfect $\mu$-th
power of an order $r$ form. About a decade ago, the second author
(JC) had constructed such a covariant using Wronskians. It will be
described below in \S\ref{definition.alphaF}; but tentatively let us denote it by $\Gott_{r,d}(F)$. Subsequently, he learnt
from the report of a colloquium lecture by Gian-Carlo
Rota~\cite{RotaLecture} that Hilbert~\cite{Hilbert} had already solved this
problem. Hilbert's construction (see~\S\ref{Hilb.construction}) is
based upon an entirely different idea; it will be denoted by $\Hilb_{r,d}(F)$. 

In fact, either of the constructions makes sense even if
$r$ does not divide $d$.  If we let $e = \gcd (r,d)$ and $d = e \, \mu$, then we have the property
\begin{equation} 
\Gott_{r,d}(F) = 0 \iff F = G^\mu \; \text{for some order $e$ form
$G$} \iff \Hilb_{r,d}(F) =0. 
\label{HG2} \end{equation}  
Both covariants turn out to be of degree $r+1$ and order $N = (r+1)(d-2)$. This, of
course, creates a strong \emph{presumption} that they might indeed be the
same (but see~\S\ref{misconception} below). This is our first result. 

\begin{Theorem} \sl There exists a nonzero rational scalar 
  $\kappa_{r,d}$ such that 
$\Gott_{r,d} = \kappa_{r,d} \, \Hilb_{r,d}$. \label{Theorem.HGeq} \end{Theorem} 
The proof will be given in \S\ref{section.Hilb.construction}. When
$r=1$, either covariant reduces to the Hessian. 

\subsection{} For $p \ge 0$, let $S_p$ denote the $(p+1)$-dimensional
space of order $p$ forms in $\ux$. We have an embedding 
\[ \P S_e \lra \P S_d,  \qquad [G] \lra [G^\mu] \] 
whose image $X = X_{e,d}$ is the variety of binary $d$-ics which are perfect
$\mu$-th powers of order $e$ forms. (In particular, $X_{1,d}$ is the rational
normal $d$-ic curve.)  Let $R = \complex[a_0, \dots, a_d]$ denote the co-ordinate ring of $\P S_d
\simeq \P^d$. Write 
\[ \Hilb_{r,d}(F) = \sum\limits_{i=0}^{N} \,
\binom{N}{i} \, h_i \, x_1^{N-i} \, x_2^i, \] 
and let $J = (h_0,h_1, \dots, h_N) \subseteq R$ denote the ideal generated
by the coefficients of $\Hilb_{r,d}$ (or what is the same, $\Gott_{r,d}$). By construction, 
the zero locus of $J$ is precisely $X$. We show that, when $r$ divides
$d$, the ideal $J$ defines $X$ as a scheme. 

\begin{Theorem} \sl 
Assume that $r$ divides $d$. Then the saturation of $J$ coincides with the defining ideal $I_X \subseteq
R$. 
\label{Theorem.Jsaturation} \end{Theorem} 
The proof will be given in \S\ref{proof.Jsaturation}. For the
case $r=1$, this theorem appears in~\cite{AluffiFaber}. 

\subsection{} In \S\ref{section.Gott.cov}, we use the plethysm
decomposition of $SL_2$-representations to exhibit $\Gott_{r,d}$ as a
special case of a family of covariants which vanish on $X$. We
baptise them as G{\"o}ttingen covariants, to commemorate the
G{\"o}ttingen school of which Hilbert was a distinguished member for
nearly five decades. In
\S\ref{GottPsi.formula}, we give an algorithm for the
symbolic computation of these covariants. The examples in 
\S\ref{section.rnc.ideal}-\S\ref{section.g26} suggest the conjecture
that the G{\"o}ttingen covariants generate all of $I_X$ when $r$
divides $d$. 

The $J$-ideals seem to obey complicated containment relations for
varying values of $r$, and there is much here that we do not understand. We give a preliminary result in this direction in
Proposition~\ref{prop.containmentJ}. The example in
\S\ref{section.twistedcubic} shows that when $r$ does not divide $d$,
the Hilbert covariants can create interesting nonreduced scheme structures on $X$. 

\subsection{} 
The problem discussed in \S\ref{power.problem} makes sense in any number of
variables. There is a classical construction due to Clebsch called the
`transfer principle',  which allows us to lift the binary solution to
$n$-ary forms. We explain this in \S\ref{section.Clebsch}, and 
construct a concomitant $\tGott_{r,d}$ of $n$-ary $d$-ics which
has exactly the same vanishing property that $\Hilb_{r,d}$ does for
binary forms (see Theorem~\ref{theorem.clebsch}). For instance, let $F$
denote a quartic form in three variables
$x_1,x_2,x_3$, which we write symbolically as 
\[ F =  a_\ux^4 = b_\ux^4 = c_\ux^4.  \] 
Then $F$ is the square of a quadratic form, if and only if the
concomitant 
\[ \tGott_{2,4} = (a \, b \, u) \, (a \, c \, u)^2 \, a_\ux \, b_\ux^3 \, c_\ux^2, \] 
vanishes on $F$. 

\subsection{} Although the Hilbert covariants were defined over a
century ago, they do not seem to have been studied much in the
subsequent years.\footnote{There is a later short note by Brioschi
  \cite{Brioschi}, but it is mostly a report on Hilbert's original
  paper and contains little that is new.} This may be partly due to Hilbert
himself, whose papers around 1890 in the \emph{Mathematische Annalen}
changed the texture of modern algebra, and to some extent caused the earlier
themes to be seen as \emph{pass{\' e}} (cf. \cite[\S II]{Fisher}). We are convinced, however, that
these covariants (and their generalisation, namely the G{\"o}ttingen
covariants) encapsulate a large amount of hitherto unexplored algebraic
geometry. 

\section{Preliminaries} 
In this section we establish notation, and explain the necessary
preliminaries in the invariant theory of binary forms. Since the
latter are less widely known now than they were a century ago, we have
included rather more background material. Some of the classical
sources for this subject are~\cite{Glenn, GY, Hilbert2, Salmon}, whereas more modern treatments
may be found in~\cite{Dolgachev1, KungRota, Olver, Processi,
Sturmfels}. In particular, for explanations pertaining to the symbolic
calculus the reader is also referred to \cite[\S 2]{Abd1}.

\subsection{$SL_2$-representations} 
The base field will be $\complex$. Let $V$ denote a two-dimensional 
complex vector space with basis $\ux = \{x_1,x_2\}$, and a natural action of 
the group $SL(V) \simeq SL_2$. For $p \ge 0$, let $S_p = \text{Sym}^p \, V$ denote the 
$(p+1)$-dimensional space of binary $p$-ics in $\ux$. Recall that $\{S_p:
p \ge 0\}$ is a complete set of finite-dimensional irreducible $SL_2$-representations,
and each finite-dimensional representation is a direct sum of
irreducibles. The reader is referred to \cite[\S 6]{FH} and
\cite[\S I.9]{Knapp} for the elementary theory of
$SL_2$-representations. 
For brevity, we will write $S_p(S_q)$ for $\text{Sym}^p (S_q)$ etc. 

\subsection{Transvectants} \label{section.trans} 
Given integers $p,q \ge 0$, we have a decomposition of representations 
\begin{equation} 
S_p \otimes S_q \simeq \bigoplus\limits_{k=0}^{\min(p,q)} \, 
S_{p+q-2k}.   \label{Clebsch-Gordan} \end{equation}
Let $A,B$ denote binary forms in $\ux$ of respective orders $p,q$. The $k$-th transvectant 
of $A$ with $B$, written $(A,B)_k$, is defined to be the image of 
$A \otimes B$ via the projection map 
\[ \pi_k: S_p \otimes S_q \lra S_{p+q-2k} \, .  \] 
It is given by the formula 
\begin{equation} (A,B)_k = \frac{(p-k)! \, (q-k)!}{p! \, q!} \, 
\sum\limits_{i=0}^k \; (-1)^i \binom{k}{i} \, 
\frac{\partial^k A}{\partial x_1^{k-i} \, \partial x_2^i} \, 
\frac{\partial^k B}{\partial x_1^i \, \partial x_2^{k-i}} \ . 
\label{trans.formula} \end{equation} 
Usually $k$ is called the index of transvection.  By convention,
$(A,B)_k = 0$, if $k > \min \, (p,q)$.
If we symbolically write $A = a_\ux^p, B = b_\ux^q$ as in \cite[Ch.~I]{GY}, then 
$(A,B)_k = (a \, b)^k \, a_\ux^{p-k} \, b_\ux^{q-k}$. A useful method for calculating transvectants of
symbolic expressions is given in \cite[\S 3.2.5]{Glenn}. 

There is a canonical isomorphism of representations 
\begin{equation} 
 S_p \stackrel{\sim}{\lra} S_p^* \, ( \, = \text{Hom}_{SL(V)}(S_p,S_0)) 
\label{self-duality} \end{equation} 
which sends $A \in S_p$ to the functional $B \lra (A,B)_p$. It is
convenient to identify each $S_p$ with its dual via this isomorphism, unless it is necessary to maintain a
distinction between them. 

\subsection{The Omega Operator} \label{section.omegaop} 
If $\ux = \{x_1,x_2\}$ and $\uy = \{y_1,y_2\}$ are two sets of binary variables, then the corresponding Omega
operator is defined to be 
\[ \Omega_{\ux \, \uy} = \frac{\partial^2}{\partial x_1 \, \partial y_2}
- \frac{\partial^2}{\partial x_2 \, \partial y_1}. \] 
Given forms $A, B$ as above, 
\[ (A,B)_k = \frac{(p-k)! \, (q-k)!}{p! \, q!} \, \left\{ \Omega_{\ux \,
  \uy}^k  \left[ A(\ux) \, B(\uy) \right] \right\}_{\uy:=\ux}. \] 
That is to say, change the $\ux$ to $\uy$ in $B$, operate $k$-times by
$\Omega$, and then revert back to the $\ux$. 
\subsection{Covariants} \label{section.covariants} 
We will revive an old notation due to Cayley, and 
write $(\alpha_0,\dots,\alpha_n \cb u,v)^n$ for the
expression
\[\sum\limits_{i=0}^n \; \binom{n}{i} \, \alpha_i \, u^{n-i} v^i. \] 
In particular 
\begin{equation} \F = (a_0,\dots,a_d \cb x_1,x_2)^d 
\label{F.gen} \end{equation} 
denotes the {\sl generic} $d$-ic, which we identify with the natural trace form in $S_d^* \, \otimes \, S_d$. 
Using the duality in~(\ref{self-duality}), this amounts to the identification 
of $a_i \in S_d^*$ with $\frac{1}{d!} \, x_2^{d-i} \, (-x_1)^i$; but
it is convenient to think of the $\ua = \{a_0, \dots, a_d\}$ as independent variables. Let $R$ denote the symmetric algebra 
\[ \bigoplus\limits_{m \ge 0} \, S_m(S_d^*) = \bigoplus\limits_{m \ge 0} \, R_m = 
\complex \, [a_0,\dots,a_d], \] 
so that $\Proj R = \P S_d \simeq \P^d$. 

Consider an $SL(V)$-equivariant embedding 
\[ S_0 \hookrightarrow R_m \otimes S_q. \] 
Let $\Phi$ denote the image of $1$ via this map, then we may write 
\begin{equation} \Phi = (\varphi_0,\dots,\varphi_q \cb x_1,x_2)^q, 
\label{Phi.cov} \end{equation} 
where each $\varphi_i$ is a homogeneous degree $m$ form in the 
$\ua$. One says that $\Phi$ is a covariant of degree $m$ and order
$q$ (of the generic $d$-ic $\F$). In other words, the space 
\[ \text{Span} \, \{\varphi_0, \dots, \varphi_q\}  \subseteq R_m \] 
is an irreducible subrepresentation isomorphic to $S_q$. 
The {\sl weight} of $\Phi$ is defined to be $\frac{1}{2}(d \, m-q)$ 
(which is always a nonnegative integer). 

In particular, $\F$ itself is a
covariant of degree $1$ and order $d$. A covariant of order $0$ is called an invariant. 
Any transvectant of two covariants is also one, hence expressions such as 
\[ (\F,\F)_4, \quad (\F,(\F,\F)_2)_3, \quad ((\F,\F)_2, (\F,\F)_4)_5, \dots  \] 
are all covariants. The Hessian coincides with $(\F,\F)_2$ up to a scalar. 
A fundamental result due to Gordan says that each covariant is a $\complex$-linear combination of such compound
transvectants (see~\cite[\S 86]{GY}). The weight of a compound
transvectant is the sum of transvection indices occurring in it; for
instance, $((\F,\F)_2, (\F,\F)_4)_5$ is of weight $2+4+5=11$. 

\subsection{} \label{ex.isobaric} 
Recall that a homogeneous form in $R$ is called isobaric
of weight $w$, if for each monomial $\prod\limits a_k^{n_k}$ appearing
in it, we have $\sum\limits_k k \, n_k = w$. If $\Phi$ is a covariant of degree-order $(m,q)$, then its coefficient
$\varphi_k$ is isobaric of weight $\frac{1}{2}(d \, m-q)+k$. For instance,
let $d=6$, and $\Phi = (\F,(\F,\F)_2)_1$, which is a covariant of degree
$3$, order $3d-6$, and hence weight $3$. Its expression begins as 
\[ \begin{aligned} 
\Phi = \; & (a_0^2 \, a_3+2 \, a_1^3-3 \, a_0 \, a_1 \, a_2) \,
x_1^{12} \,  + \\ 
& (12 \, a_1^2 \, a_2-15 \, a_0 \, a_2^2+3 \, a_0^2 \, a_4) \, x_1^{11} \,
x_2 \, + \\ &
(15 \, a_1 \, a_2^2+3 \, a_0^2 \, a_5+18 \, a_0 \, a_1 \, a_4+24 \,
a_1^2 \, a_3-60 \, a_0 \, a_2 \, a_3) \, x_1^{10} \, x_2^2  \, + \\
& ( 25 \, a_2^3+60 \, a_1^2 \, a_4-80 \, a_0 \, a_3^2+a_0^2 \, a_6-30 \,
a_4 \, a_0 \, a_2+24 \, a_1 \, a_0 \, a_5) \, x_1^9 \, x_2^3  \, + \dots, 
\end{aligned} \] 
and one sees that the successive coefficients are isobaric of weights
$3,4,5$ etc. 

\subsection{The Cayley-Sylvester formula} Let $C(d, m,q)$ denote the
vector space of covariants of degree-order $(m,q)$ for binary $d$-ics;
its dimension is the same
as the multiplicity of $S_q$ in the irreducible decomposition of $R_m
\simeq S_m(S_d)$. This number is given by the Cayley-Sylvester formula (see~\cite[Corollary
4.2.8]{Sturmfels}). For integers $n,k,l$, let $\pi(n, k, l)$ denote
the number of partitions of $n$ into $k$ parts such that no part
exceeds $l$. Then 
\[ \zeta(d,m,q) = \dim C(d,m,q) = \pi \left(\frac{d \, m-q}{2},d,m \right) -
\pi \left(\frac{d \, m-q-2}{2},d,m \right). \] 
For instance, $\zeta(6,3,6)=\pi(6,6,3)-\pi(5,6,3)=7-5=2$, and it is easy to check (e.g., by
specialising $\F$) that 
\[ \F \, (\F,\F)_6, \quad (\F, (\F,\F)_4)_2, \] is a basis of $C(6,3,6)$. 

\subsection{} \label{misconception} 
This is perhaps the correct place to forestall one possible misconception about
Theorem~\ref{Theorem.HGeq}. Recall that $\Gott_{r,d}$ and $\Hilb_{r,d}$ both
have degree $r+1$ and order $N = (r+1) \, (d-2)$. If it were the case that
\begin{equation} 
\zeta(d,r+1,N)=1, 
\label{zetaeq1} \end{equation} 
then one could immediately conclude that the two 
must be equal up to a scalar. But such may not be the case. For instance, if
$r=5,d=15$, then $\zeta(15,6,78) = 4$. 
Hence, Theorem~\ref{Theorem.HGeq} does not follow from general
multiplicity considerations, but instead requires an explicit hard
calculation. However, (\ref{zetaeq1}) is true for $r=1,2$. (This
can be seen from the plethysm formulae in~\cite[\S I.8]{MacDonald}.) 

\subsection{} An alternate equivalent definition of a covariant is as follows. Let
$g = \left( \begin{array}{cc} \alpha & \gamma \\ \beta &
    \delta \end{array} \right)$, where $\alpha, \dots, \delta$ are
regarded as independent indeterminates. Write 
\[ x_1 = \alpha \, x_1' + \beta \, x_2', \quad x_2 = \gamma \, x_1' +
\delta \, x_2',  \] 
and substitute into (\ref{F.gen}). Determine expressions $a_0',\dots,
a_d'$ such that we have an equality 
\[ (a_0',\dots,a_d' \cb x_1',x_2')^d = (a_0,\dots,a_d \cb
x_1,x_2)^d;  \] 
then each $a_i'$ is a polynomial expression in the $\ua$ and $\alpha,
\dots, \delta$. Now let $\Phi \in \complex[a_0, \dots, a_d; x_1,x_2]$ be a bihomogeneous form
of degrees $m,q$ respectively in $\ua, \ux$. Then $\Phi$ is a
covariant, if and only if the following identity holds: 
\begin{equation} 
\Phi(a_0', \dots, a_d'; x_1', x_2') = (\alpha \, \delta - \beta \,
\gamma)^{\frac{dm-q}{2}} \, \Phi(a_0, \dots, a_d; x_1,
x_2). \label{cov.identity} \end{equation} 

\subsection{Covariants and Differential
  Operators} \label{section.cov.diff} 
Consider the following differential operators: 
\begin{equation} 
E_{+} = \sum\limits_{i=0}^{d-1} \, (d-i) \, a_{i+1} \,
\frac{\partial}{\partial a_i},  \qquad 
E_{-} = \sum\limits_{i=1}^d i \, a_{i-1} \, \frac{\partial}{\partial
 a_i}, \qquad 
E_0 = \sum\limits_{i=0}^d \, (2i-d)  \, a_i \,
\frac{\partial}{\partial a_i}, 
\label{cayley.equations2} \end{equation} 
and 

\[ 
\Gamma_{+} = E_{+} -  x_1 \, \frac{\partial}{\partial x_2}, \qquad 
\Gamma_{-} = E_{-} - x_2 \, \frac{\partial}{\partial x_1}, \qquad 
\Gamma_0 = E_0 + (x_1 \, \frac{\partial }{\partial x_1} - 
x_2 \, \frac{\partial }{\partial x_2}). \] 
\begin{Proposition} \sl 
A bihomogeneous form $\Phi$ is a covariant, if and only if 
\begin{equation} 
\Gamma_{+} \, \Phi = \Gamma_{-} \, \Phi = \Gamma_{0} \, \Phi =
0. \label{GammaPhi} \end{equation} 
\end{Proposition} 
A proof is given in~\cite[\S 149]{Salmon} (also see~\cite[\S
4.5]{Sturmfels}), but the central idea is the following: $\Phi$ is a covariant exactly when it
remains unchanged by an $SL_2$-action, i.e., when it is annihilated
by the Lie algebra ${\mathfrak{sl}}_2$. Let 
\[ J_{+} = \left( \begin{array}{rr} 0 & 1 \\ 0 & 0 \end{array} \right), \quad 
J_{-} = \left( \begin{array}{rr} 0 & 0 \\ 1 & 0 \end{array} \right), \quad 
J_0 = \left( \begin{array}{rr} 1 & 0 \\ 0 & -1 \end{array} \right), \] 
denote the standard generators of ${\mathfrak{sl}}_2$. Choose a path 
$t \lra g_t$ starting from the identity element in $SL_2$,
and apply condition (\ref{cov.identity}) to $g_t$. For the three
cases $J_\star = \left[\frac{d g_t}{d t}\right]_{t=0}$ where $\star \in
\{+,-,0\}$, we respectively get the identities in (\ref{GammaPhi}). \qed 

\smallskip 

The first coefficient $\varphi_0$ is called the source (or
seminvariant) of $\Phi$. From~(\ref{GammaPhi}), we get equations 
\begin{equation} 
E_{-}(\varphi_0) =0, \quad \text{and} \quad 
\varphi_k = \frac{(q-k)!}{q !}  \, E_{+}^k(\varphi_0) \quad \text{for $0 \le k \le
 q$.} 
\label{source.eqns} \end{equation} 
Thus one can recover the entire covariant from the source alone. 
Moreover, a homogeneous isobaric form $\psi$ in the $\ua$ can be a
source (of some covariant), if and only if it satisfies the condition $E_{-}(\psi)=0$. 

The commutation relations between the $E_\star$ are parallel to the ones between the standard generators of
${\mathfrak{sl}}_2$, i.e., 
\[ [E_{+}, E_{-}] = E_0, \quad [E_0, E_{+}] = 2 \, E_{+}, \quad [E_0,
E_{-}] = - 2 \, E_{-}.  \] 
The following lemma will be needed in \S\ref{Hilb.construction}. 
\begin{Lemma} \sl 
For $n \ge 0$, we have an identity 
\[ E_{-} \, E_{+}^{n+1} = E_{+}^{n+1} \, E_{-} - 
(n+1) \, E_{+}^n \, E_0 - n \, (n+1) \, E_{+}^n. \] 
\label{lemma.E} \end{Lemma} 
\demo This follows by a straightforward induction on $n$. \qed 

\subsection{Wronskians} 
Let $m,n \ge 0$ be integers such that $m \le n+1$. 
Consider the following composite morphism of representations 
\[ w: \wedge^m S_n \stackrel{\sim}{\lra} S_m(S_{n-m+1}) \lra S_{m(n-m+1)}, \] 
where the first map is an isomorphism (described in~\cite[\S2.5]{AC2}) and the second is the natural surjection. 
Given a sequence of binary $n$-ics $A_1,\dots,A_m$, define 
their Wronskian $W(A_1,\dots,A_m)$ to be the image $w(A_1 \wedge \dots
\wedge A_m)$. It is given by the determinant 
\[ (i,j) \lra \frac{\partial^{m-1}  \, A_i}{\partial x_1^{m-j} \, \partial \, x_2^{j-1}}, \quad 
(1 \le i,j \le m). \] 
The $\{A_i\}$ are linearly dependent over $\complex$, if and only if
$W(A_1, \dots, A_m)=0$. (The `only if' part is obvious. For 
the converse, see~\cite[\S 1.1]{Meulien}.) 

\section{The G{\"o}ttingen covariants} \label{section.Gott.cov} 

\subsection{} \label{section.defn.rde}
Henceforth assume that $r,d$ are positive integers, and let $e =
\gcd(r,d)$. Write $d = e \, \mu$ and $r = e \, \mu'$. Consider the embedding 
\[ \P S_e \stackrel{\imath}{\lra} \P S_d, \quad [G] \lra [G^\mu]. \] 
Let $X_{e,d}$ denote the image variety. We have a factorisation 
\[ \diagram 
{} & \; \P S_\mu(S_e) \ar@{.>}[d]^\pi \\ 
\P S_e \ar[ur]^{v_\mu} \ar[r]^\imath & \P S_d 
\enddiagram 
\] 
where $v_\mu$ is the $\mu$-fold Veronese embedding, and $\pi$ is the
projection coming from the surjective map $S_\mu(S_e) \lra S_{e \, \mu} =
S_d$. Thus $\imath$ corresponds to the incomplete linear series $S_d
\subseteq H^0({\mathcal O}_{\P S_e}(\mu))$.

\subsection{} \label{definition.alphaF} 
In this section we will define the covariants
$\Gott_{r,d}$. For $F \in S_d$, we have a morphism 
\[ \alpha_F : S_r \lra S_{r+d-2}, \quad A \lra (A,F)_1 = \frac{1}{rd}
\, \left| \begin{array}{cc} A_{x_1} & A_{x_2} \\ F_{x_1} &
  F_{x_2} \end{array} \right|,  \] 
where $A_{x_i}$ stands for $\frac{\partial A}{\partial x_i}$ etc. 

\begin{Proposition} \sl With notation as above, 
\[ \ker \alpha_F \neq 0 \iff [F] \in X_{e,d}. \] 
\end{Proposition} 
\demo Assume $F = G^\mu$ for some $G$, then 
$\alpha_F(G^{\mu'}) = (G^{\mu'}, G^\mu)_1 =0$. 

Alternately, assume that $(A,F)_1=0$ for some nonzero $A$. We will
construct a form $G$ such that (up to scalars) $A = G^{\mu'}, F =
G^\mu$. Let $\ell \in S_1$ be a linear form which divides
either $A$ or $F$; after a change of variables we may assume $\ell =
x_1$. Suppose that $a,f$ are the highest powers of $x_1$ which
divide $A,F$ respectively, and write $A = x_1^a \, \tA, F = x_1^f
\, \tF$. Starting from the relation
$A_{x_1} \, F_{x_2} = A_{x_2} \, F_{x_1}$, after expanding and rearranging the terms, we get 
\[ x_1^{a+f} \, (\tA_{x_1} \, \tF_{x_2} - \tA_{x_2} \, \tF_{x_1}) =
x_1^{a+f-1} \, (-a \, \tA \, \tF_{x_2} + f \, \tA_{x_2} \, \tF),  \] 
hence $x_1$ must divide $(-a \, \tA \, \tF_{x_2} + f \,
\tA_{x_2} \, \tF)$. Thus, either $\tA_{x_2} = \tF_{x_2} =0$ (and so
$a=r,f=d$), or the terms with highest powers of $x_2$ in $a \,
\tA \, \tF_{x_2}$ and $f \, \tA_{x_2} \, \tF$ cancel against each
other. In the latter case, 
\[ a \, (d-f) = f \, (r-a) \implies a \, d = f \, r \implies a \, \mu
= f \, \mu'. \] 
In either case, $\mu' \, | \,  a$ and $\mu \, | \, f$. Define $G$
such that $x_1$ appears in it exactly to the power $\frac{f}{\mu}$, and similarly
for all $\ell$. \qed 

\medskip

Now consider the composite morphism 
\[ 
S_0 \simeq \wedge^{r+1} S_r \stackrel{\wedge^{r+1} \alpha_F}{\lra} 
\wedge^{r+1} S_{r+d-2} \simeq S_{r+1}(S_{d-2}) \lra S_{(r+1)(d-2)}. 
\] 
The image of $1 \in S_0$ is the Wronskian $W(\alpha_F(x_1^r),
\alpha_F(x_1^{r-1} x_2), \dots, \alpha_F(x_2^r))$,  
which we define to be $\Gott_{r,d}(F)$. To recapitulate, 
$\Gott_{r,d}(F)$ is the determinant of the $(r+1) \times (r+1)$ matrix 
\begin{equation} 
(i,j) \lra \frac{\partial^r \, C_i}{\partial x_1^{r-j} \, \partial
  \, x_2^j}, \quad 
(0 \le i,j \le r). \label{formula.Grd} \end{equation} 
where $C_i = (x_1^{r-i} \, x_2^i,F)_1$. Then 
\[ \Gott_{r,d}(F) = 0 \iff \ker \alpha_F \neq 0 \iff [F] \in X_{e,d}. \]  
Each matrix entry is linear in the coefficients of $F$,
and of order $d-2$ in $\ux$, hence $\Gott_{r,d}$ has degree $r+1$ and order
$N = (r+1)(d-2)$. 

In the next section, we will generalise this construction to obtain a
family of covariants vanishing on $X_{e,d}$. The reader who is
more interested in Hilbert's solution may proceed directly to
\S\ref{Hilb.construction}. 

\subsection{}
Let 
\[ \B = (b_0, b_1, \dots, b_{d-2} \cb x_1,x_2)^{d-2}, \] 
denote a generic form of order $d-2$, with a new set of
indeterminates $\ub$.  As in \S\ref{section.covariants}, the
$\ub$ can be seen as forming a basis of $S_{d-2}^* \simeq
S_{d-2}$. Let $\Psi (\ub, \ux)$ denote a covariant of degree $r+1$
and order $q$ of $\B$. Then $\Psi$ corresponds to an embedding $S_q \lra
S_{r+1}(S_{d-2})$, which can be described as follows: if we realise
$S_{r+1}(S_{d-2})$ as the space of degree $r+1$ forms in the $\ub$, then $A(\ux) \in S_q$ gets sent to
$(A,\Psi)_q$. After dualising, we get a morphism 
\[ f_\Psi: S_{r+1}(S_{d-2}) \lra S_q. \] 
Now consider the composite morphism 
\[ S_0 \simeq \wedge^{r+1} S_r \stackrel{\wedge^{r+1} \alpha_F}{\lra}
\wedge^{r+1} S_{r+d-2} \simeq S_{r+1}(S_{d-2})
\stackrel{f_\Psi}{\lra} S_q,  \] 
and let $\Gott_\Psi(\F)$ denote the image of $1 \in S_0$, which will
be called the G{\"o}ttingen covariant of $\F$ associated to $\Psi$. It is of
the same degree and order as $\Psi$, and hence its weight is $(r+1)$ more than that of $\Psi$.  
In particular, $\Gott_{r,d}(\F)$ is the same as $\Gott_{\B^{r+1}}(\F)$. As before, 
\begin{equation}  [F] \in X_{e,d} \implies \Gott_\Psi(F) = 0.  \label{Psi.implications} \end{equation} 

\subsection{The calculation of $\Gott_\Psi$} \label{GottPsi.formula}
One can calculate $\Gott_\Psi(\F)$ explicitly by following the sequence of
maps above, which amounts to the following recipe: 
\begin{itemize} 
\item Introduce $2 \, (r+1)$ sets of binary variables 
\[ \uy_{(i)} = \{y_{i1},  y_{i2} \}, \quad \uz_{(i)} =
\{z_{i1},z_{i2}\} \quad \text{for  \; $0  \le i 
  \le r$; } \] 
and let $\Omega_{\uy_{(i)} \, \uz_{(i)}}$ be the corresponding Omega
operators. 
\item 
Let $\cW$ denote the determinant 
\[ (i,j) \lra \left. \frac{\partial^r (x_1^{r-i} \, x_2^i,\F)_1}{\partial
  x_1^{r-j} \, \partial x_2^j} \right|_{\ux \lra \uy_{(i)}}, \quad (0 \le i, j \le r). \] 
This is similar to~(\ref{formula.Grd}), except that the $\uy_{(i)}$ variables are
used throughout the $i$-th row. Let $\cW^\sharp$ denote the
symmetrisation of $\cW$ with respect to the sets $\uy_{(i)}$, i.e., 
\[ \cW^\sharp = \sum\limits_\sigma \; \cW(\uy_{\sigma(0)}, \dots,
\uy_{\sigma(r)}), \] 
the sum quantified over all permutations $\sigma$ of $\{0,
\dots, r\}$. Then $\cW^\sharp$ is of degree
$r+1$ in $\ua$, and of order $d-2$ in each $\uy_{(i)}$. 
\item Write 
\[ \Psi = (\psi_0,\dots, \psi_q \cb x_1,x_2)^q, \] 
where each $\psi_i$ is a degree $r+1$ form in $\ub = \{b_0, \dots,
b_{d-2} \}$. 
Introduce $r+1$ sets of variables $\ub_{(0)}, \dots, \ub_{(r)}$, where 
\[ \ub_{(i)} = \{b_{i \, 0}, \dots, b_{i \, d-2} \}, \] 
and let $\widetilde \Psi$ be the total polarisation of $\Psi$ with respect to
the new variables (see~\cite[\S 1.1]{Dolgachev2}). Then $\widetilde \Psi$ is linear in each set
$\ub_{(i)}$. 
\item Let $\widehat \Psi$ denote the form obtained from $\widetilde \Psi$ by replacing
  $b_{ik}$ with $\frac{1}{(d-2)!} \, z_{i2}^{d-2-k}(-z_{i1})^k$ for $0 \le i \le r$ and
  $0 \le k \le d-2$. (This is similar to the identification of $a_i$
  as in \S\ref{section.covariants}.) Thus $\widehat \Psi$ is of order $q$ in $\ux$,
  and of order $d-2$ in each $\uz_{(i)}$. 
\item 
Finally, 
\[ \Gott_\Psi(\F) = [ \, \Omega_{\uy_{(0)} \, \uz_{(0)}}^{d-2} \circ \dots \circ
\Omega_{\uy_{(r)} \, \uz_{(r)}}^{d-2} \, ] \, \widehat\Psi \, \cW^\sharp \,.  \] 
This removes all the $\uy_{(i)}$ and $\uz_{(i)}$ variables,  which leaves a
form of degree $r+1$ in the $\underline{a}$ and order $q$ in $\ux$. 
\end{itemize}

\subsection{} It may happen that $\Gott_\Psi$ is identically zero, 
even if $\Psi$ is nontrivial. (Hence the implication
in~(\ref{Psi.implications}) is not reversible in general.) For
instance, recall that a generic binary $d$-ic has a cubic 
invariant exactly when $d$ is a multiple of $4$. Now let
$r=2$, and assume $d \equiv 2 \; (\text{mod 4})$. Then $\Psi
=(\B,(\B,\B)_{\frac{d-2}{2}})_{d-2}$ is a nontrivial cubic invariant of
$\B$, but $\Gott_\Psi$ must vanish identically. 

If $d$ is a divisor of $r$, then $X_{e,d} = \P S_d$, and
in that case all $\Gott_\Psi$ are identically zero. 

\medskip 

\noindent {\bf N.B.} Henceforth, if $A,B$ are two quantities, we will write $A \doteq B$ to
mean that $A =  c \, B$ for some unspecified nonzero rational scalar
$c$. This will be convenient in symbolic calculations, where more and
more unwieldy scalars tend to accumulate at each stage. 

\subsection{} \label{Gott.quadratic} 
As an example, we will follow this recipe when 
$r=1$ and $\Psi$ is any quadratic covariant. Write symbolically
\[ \F = \alpha_\ux^d = \beta_\ux^d, \qquad \B = p_\ux^{d-2} = q_\ux^{d-2}. \] 
Every quadratic covariant of $\B$ must be of the form 
\[ \Psi = (\B,\B)_{2n} = (p \, q)^{2n} \, p_\ux^{d-2-2n} \,
q_\ux^{d-2-2n}, \] 
for some $n$ in the range $0 \le n  \le \frac{d-2}{2}$. 
Using $\alpha, \beta$ for the two rows of $\cW$, we get 
\[ \cW \doteq \left| \begin{array}{rr} 
\alpha_1 \, \alpha_2 \, \alpha_{\uy_{(0)}}^{d-2} & \alpha_2^2 \, \alpha_{\uy_{(0)}}^{d-2}
\\ 
\beta_1^2 \, \beta_{\uy_{(1)}}^{d-2} & \beta_1 \, \beta_2 \,
\beta_{\uy_{(1)}}^{d-2} \end{array} \right| = \alpha_2 \, \beta_1 \, (\alpha
\, \beta) \, \alpha_{\uy_{(0)}}^{d-2} \, \beta_{\uy_{(1)}}^{d-2}, \] 
and hence 
\[ \cW^\sharp \doteq (\alpha \, \beta)^2 \, \alpha_{\uy_{(0)}}^{d-2} \,
\beta_{\uy_{(1)}}^{d-2}. \] 
Now the symbolic expression for $\widetilde \Psi$ is the same as the
one for $\Psi$, once we make the convention that $p,q$ respectively refer
to the $\ub_{(0)}, \ub_{(1)}$ variables. Then 
\[ \widehat \Psi \doteq (\uz_{(0)} \, \uz_{(1)})^{2n} \, (\ux \, \uz_{(0)})^{d-2-2n} \, (\ux
\, \uz_{(1)})^{d-2-2n};  \] 
and finally, 
\[ \Gott_\Psi \doteq (\alpha \, \beta)^{2n+2} \, \alpha_\ux^{d-2-2n} \, \beta_\ux^{d-2-2n}. \] 
We have proved the following: 
\begin{Proposition} \sl 
If $\Psi = (\B,\B)_{2n}$, then $\Gott_{\Psi} \doteq (\F, \F)_{2n+2}$. 
\label{prop.quadG} \end{Proposition} 

In particular, $\Gott_{1,d} \doteq \Gott_{\B^2} = 
\Gott_{(\B,\B)_0} \doteq (\F, \F)_2$ is the
Hessian of $\F$. Similar calculations show that 
\begin{equation} 
\begin{array}{l} 
\Gott_{2,d} \doteq (\F, (\F,\F)_2)_1, \\ 
\Gott_{3,d} \doteq 3 \, (2 \, d-3) \, (\F,\F)_2^2 - 2 \, (d-2) \, \F^2 \,
(\F,\F)_4, \\ 
\Gott_{4,d} \doteq 2 \, (3 \, d-4) \, (\F,\F)_2 \, (\F, (\F,\F)_2)_1 - (d-3) \,
\F^2 \, (\F,(\F,\F)_4)_1.  \end{array} 
\label{Gott.lowr} \end{equation} 
(Such formulae are derived for the $\Hilb_{r,d}$ in~\cite{Brioschi}
and~\cite{Hilbert}, but this makes no difference in view of
Theorem~\ref{Theorem.HGeq}.) However, as $r$ grows, it quickly
begins to get more and more tedious to execute this recipe. 

\section{Hilbert's construction} 
\label{section.Hilb.construction} 
In this section we will describe Hilbert's construction of his
covariants $\Hilb_{r,d}$, and later prove that the outcome
coincides with $\Gott_{r,d}$ up to a scalar. 

The underlying idea is as follows. Suppose, for instance, that 
$F$ is an order $10$ form such that $F = G^5$ for some quadratic $G$. Then substituting $x_1=1, x_2=z$, we have 
\[ \frac{d^3}{d z^3} \;  \sqrt[5]{F(1,z)} = 0. \] 
One should like to convert the left-hand side into a covariant
condition on $F$; but this requires some technical modifications. We begin
by constructing the source of Hilbert's covariant. 

\subsection{} \label{Hilb.construction} 
Define 
\[ h_0 = a_0^{r+1-\frac{r}{d}} \, E_+^{\, r+1} \,
(a_0^{\, \frac{r}{d}}).   \] 
This is easily seen to be an isobaric homogeneous form of degree and weight
$r+1$ in the $\ua$. For instance, 
\[ h_0 = 
\begin{cases} (d-1) \, (a_0 \, a_2 - a_1^2) & \text{if $r=1$,} \\ 
(2 d^2 -6d+4) \, a_0^2 \, a_3 -(6d^2-18d+12) \, a_0 \, a_1 \, a_2 +
(4 \, d^2- 12 \, d + 8) \, a_1^3 & \text{if $r=2$.} 
\end{cases} \] 
\begin{Lemma} \sl 
The form $h_0$ is a source. 
\end{Lemma} 
\demo We want to show that $E_{-} \, h_0 = a_0^{r+1-\frac{r}{d}} \,
E_{-} \, E_{+}^{\, r+1}(a_0^{\, \frac{r}{d}})$ vanishes. Apply Lemma~\ref{lemma.E}, and
note that 
\[ E_0 \, (a_0^{\frac{r}{d}})= - r \, a_0^{\frac{r}{d}}, \quad E_{-} (a_0^{\frac{r}{d}})=0, \] 
which implies the result. \qed 

\smallskip 

Since $h_0$ has weight $r+1$, the covariant corresponding to $h_0$
must have order $N = (r+1)(d-2)$. The Hilbert covariant is defined to be 
\begin{equation} 
\Hilb_{r,d}(\F) = (h_0, \dots, h_N \cb x_1,
x_2)^N, \label{Hilb.cov.definition} \end{equation} 
where 
\begin{equation} 
h_k = \frac{(N-k)!}{N !}  \, E_{+}^k(h_0) 
\qquad \text{for $0 \le k \le N$.} 
\label{formula.hk} \end{equation}  

\subsection{} In order to prove Theorem~\ref{Theorem.HGeq}, it will
suffice to show that $\Hilb_{r,d}$ and $\Gott_{r,d}$ have the same
source up to a scalar. We will avoid writing such scalars explicitly in the course of the
calculation, but see formula~(\ref{formula.krd}) below. 

Let $\ua = (a_0, \dots, a_d)$ denote a $(d+1)$-tuple of
complex variables. For $t \in \complex$ define 
\[ \gamma_t: \complex^{d+1} \lra \complex^{d+1}, \quad 
(a_0,a_1, \dots, a_d) \lra (a_0(t), a_1(t), \dots, a_d(t)),  \] 
by the formula 
\[ (a_0, \dots, a_d \cb 1,z+t)^d = (a_0(t), \dots, a_d(t) \cb
1,z)^d. \] 
It is easy to see that 
\[ a_i(t) = a_i + (d-i) \, a_{i+1} \, t + O(t^2) \quad \text{for $0
  \le i \le d-1$}. \] 
Hence, given an analytic function $\phi: \complex^{d+1} \lra
\complex$, we have an equality 
\[ E_+ \, \phi = \left[\frac{d}{d \, t} (\phi(\gamma_t))
\right]_{t=0}. \] 
Iterating this formula, 
\begin{equation} 
E_+^n \, \phi =\left[\frac{\partial^n \phi(\gamma_{t_1+\dots +t_n})}{\partial t_1\cdots\partial t_n}\right]_{t_1=\cdots=t_n=0}
= \left[\frac{d^n \, \phi(\gamma_t)}{d \, t^n}\right]_{t=0}. 
\label{formula.Ediff} \end{equation} 
Now write $f(z) = (a_0, \dots, a_d \cb 1,z)^d$, and apply this to the function 
\[ \phi(\ua) = a_0^{\, \frac{r}{d}} = f(0)^{\frac{r}{d}} ,   \] 
which gives the expression 
\begin{equation}
h_0=f(0)^{r+1-\frac{r}{d}} \, \left[ \frac{d^{\, r+1}}{d \, t^{r+1}}
  \, f(t)^{\frac{r}{d}} \right]_{t=0},  \label{h0.formula} \end{equation} 
for the source of $\Hilb_{r,d}$. 

\subsection{} 
We make a small digression to prove that $\Hilb_{r,d}$ has the
required vanishing property. 
\begin{Proposition} \sl 
Let $F$ be a $d$-ic. Then 
\[ \Hilb_{r,d}(F) =0 \iff [F] \in X_{e,d}. \] 
\end{Proposition} 
\demo 
This will of course follow from Theorem~\ref{Theorem.HGeq}, but even so, 
we include an independent proof. After a change of variables, we may
assume $a_0 \neq 0$. Write 
\begin{equation} f(t)^{\frac{r}{d}} = a_0^{\frac{r}{d}} \; (1 + \sum\limits_{m \ge 1}
  \; \frac{\theta^m}{m!} \, t^m). 
\label{powseries.f} \end{equation}
Using the reformulation of $E_+$ above, we have 
\[ a_0^{\frac{r}{d}} \, \theta_{m+1} = E_+(a_0^{\frac{r}{d}} \,
\theta_{m}). \] 
Now a simple induction shows that there are identities 
\[ a_0^{r+1} \, \theta_{r+1} = h_0, \quad 
a_0^{r+2} \, \theta_{r+2} \doteq a_0 \, h_1 + \Box_1 \, h_0, \] 
and in general 
\[ a_0^{r+1+k} \, \theta_{r+1+k} \doteq a_0^k \, h_k + \sum\limits_{i=1}^k \; \Box_i \, h_{k-i} , \] 
for some homogeneous polynomials $\Box_i(\ua)$ of degree $k$ and
weight $i$. (Here we have set $h_i = 0$ for $i > N$.) If $h_0, h_1,
\dots $ etc.~all vanish, then so do $\theta_k$ for $k
\ge r+1$, and the power series in~(\ref{powseries.f}) becomes a polynomial of degree
$\le r$. Thus $f(t)$ reduces to a perfect $\mu$-th power. Conversely, if
$f(t) = g(t)^\mu$, then $f^{\frac{r}{d}} = g^{\mu'}$ is of degree $\le
r$, and hence $h_0 = h_1 = \dots =0$. \qed 

\subsection{} 
We should like to calculate the source of $\Gott_{r,d}$ as defined by
the determinant in~(\ref{formula.Grd}). However, the dehomogenisation
in the previous section is with respect to the other variable, so a
preparatory step is needed. The Wronskian construction is equivariant,
hence in the notation of \S\ref{section.intro1}, we have an identity 
\[ W(C_0,\dots, C_r)(x_1,x_2) = W(C_0^g, \dots, C_r^g)(x_2,-x_1), \] 
for $g = \left( \begin{array}{rr} 0 & -1
    \\ 1 & 0 \end{array} \right)$. Let $z = -x_2/x_1$, then (up to a
scalar) the right-hand side becomes 
\[(-x_1)^N \times \left|
\begin{array}{ccc}
v_0^{(r)} & \cdots & v_0\\
\vdots & & \vdots\\
v_r^{(r)} & \cdots & v_r \end{array} \right|, \]
where
\[ v_i(z)=C_i^g(z,1)=C_i(-1, z), \quad \text{and} \quad v_i^{(k)} =
\frac{d^k}{dz^k} v_i.\]
Substituting $x_1=1,x_2=0$, 
\[ g_0 = \text{source of $\Gott_{r,d}$} \doteq \left|
\begin{array}{ccc}
v_0^{(r)}(0) & \cdots & v_0(0)\\
\vdots & & \vdots\\
v_r^{(r)}(0) & \cdots & v_r(0) \end{array} \right|. \] 
By definition, 
\[ C_i \doteq \frac{\partial \, [x_1^{r-i} \, x_2^i]}{\partial x_1} 
    \frac{\partial F}{\partial x_2} - 
\frac{\partial \, [x_1^{r-i} \, x_2^i]}{\partial x_2} \, \frac{\partial
    F}{\partial x_1}. \] 
Using $F(x_1,x_2)=x_1^d \, f(\frac{x_2}{x_1})$, this can be rewritten as 
\[ C_i \doteq \left(
i \, d \, x_1^{d+r-i-1} \, x_2^{i-1} \, f\left(\frac{x_2}{x_1}\right)
-r \, x_1^{d+r-i-2} \, x_2^i \, f'\left(\frac{x_2}{x_1}\right)
\right), \]
and therefore
\[ v_i(z) \doteq \left( i \, d \, (-1)^{d+r-i-1} \, z^{i-1} \, f(-z)
  -r \, (-1)^{d+r-i-2} \, z^i \, f'(-z) \right). \] 
Let
\begin{equation} b_i(z)= i \, d \, z^{i-1} \, f(z)-r \, z^i \, f'(z), \label{definition.biz} \end{equation} 
so that $v_i(-z)=(-1)^{d+r} \, b_i(z)$. After reordering the columns, 
\[
g_0 \doteq \left|
\begin{array}{ccc}
b_0(0) & \cdots & b_0^{(r)}(0)\\
\vdots & & \vdots\\
b_r(0) & \cdots & b_r^{(r)}(0)
\end{array} \right|. \] 
Now the key step is to write 
\[ b_i \doteq f(z)^{1+\frac{r}{d}} \, \underbrace{\frac{d}{dz} \, \left(z^i
  f(z)^{-\frac{r}{d}}\right)}_{c_i}, \]
which is similar to the idea of an integrating factor in the theory of
ordinary differential equations. Then 
\begin{equation} g_0 \doteq (a_0^{1+\frac{r}{d}})^{r+1}\times \left|
\begin{array}{ccc}
c_0(0) & \cdots & c_0^{(r)}(0)\\
\vdots & & \vdots\\
c_r(0) & \cdots & c_r^{(r)}(0)
\end{array}
\right|\ .
\label{det.ci} \end{equation} 
Now let $\mu_i =z^i$ and $\omega=f(z)^{-\frac{r}{d}}$, so that
$c_i=(\mu_i \, \omega)'$. By the Leibniz rule, 
\[ c_i^{(j)} = (\mu_i \, \omega)^{(j+1)}=\sum_{k=0}^{r+1} \binom{j+1}{k} \, 
\mu_i^{(k)} \omega^{(j+1-k)} \quad \text{for $0\le j\le r$,} \] 
with the convention that $\binom{j+1}{k}=0$ if $k > j+1$. Since
$\mu_i^{(r+1)}=0$, we can stop the summation at $k=r$, and factor the determinant
in~(\ref{det.ci}) as
\begin{equation} \left|
\begin{array}{ccc}
\mu_0(0) & \cdots & \mu_0^{(r)}(0)\\
\vdots & & \vdots\\
\mu_r(0) & \cdots & \mu_r^{(r)}(0)
\end{array}
\right|\times \Delta, \label{det.mu.delta} \end{equation} 
where
\[ \Delta=
\left| \begin{array}{ccccc}
\left(\begin{array}{c} 1\\ 0\end{array}\right)\omega' &
\left(\begin{array}{c} 2\\ 0\end{array}\right)\omega'' &
\cdots & \cdots &
\left(\begin{array}{c} r+1\\ 0\end{array}\right)\omega^{(r+1)}\\
\left(\begin{array}{c} 1\\ 1\end{array}\right)\omega &
\left(\begin{array}{c} 2\\ 1\end{array}\right)\omega' &
\cdots & \cdots &
\left(\begin{array}{c} r+1\\ 1\end{array}\right)\omega^{(r)}\\
0 &
\left(\begin{array}{c} 2\\ 2\end{array}\right)\om &
\cdots & \cdots &
\left(\begin{array}{c} r+1\\ 2\end{array}\right)\om^{(r-1)}\\
\vdots & \ddots & \ddots & & \vdots\\
0 & \cdots & 0 &
\left(\begin{array}{c} r\\ r\end{array}\right)\om &
\left(\begin{array}{c} r+1\\ r\end{array}\right)\om'
\end{array}
\right|
\]
evaluated at $z=0$. The first determinant in~(\ref{det.mu.delta}) is a
pure rational constant, so it only remains to calculate $\Delta$. 
 
\subsection{} 
Let $\nu=\omega^{-1}=f(z)^{\frac{r}{d}}$. For $f(0)\neq 0$,
these are holomorphic functions of $z$ near the origin. From 
\[ \omega \, \nu = (\omega_0+\omega_1 \, z+\omega_2 \, z^2+\cdots) \,
(\nu_0+\nu_1 \, z+\nu_2 \, z^2+\cdots) =1, \]
we get the linear system 
\[
\omega_0 \, \nu_0=1, \quad \text{and} \quad 
\sum\limits_{i=0}^k \, \omega_i \, \nu_{k-i} =0, \quad 
\text{for $1 \le k \le r+1$.} \] 

Solving for $\nu_{r+1}$ by Cramer's rule, 
\begin{equation} \nu_{r+1}=\frac{1}{\omega_0^{r+2}}\times
\left|
\begin{array}{ccccc}
\omega_0 & 0 & \cdots & \cdots & 1\\
\omega_1 & \om_0 & 0 & \cdots & 0\\
\vdots &  & \ddots & & \vdots\\
\omega_r & & & \om_0 & 0\\
\omega_{r+1} & \cdots & \cdots & \omega_1 & 0
\end{array}
\right| = 
\frac{(-1)^{r+1}}{\omega_0^{r+2}} \times
\left|
\begin{array}{ccccc}
\omega_1 & \om_2 & \cdots & \cdots & \omega_{r+1}\\
\omega_0 & \omega_1 & \cdots & \cdots & \omega_{r}\\
0 & \omega_0 & \ddots & & \vdots\\
\vdots & \ddots & \ddots & \ddots & \vdots\\
0 & \cdots & 0 & \omega_0 & \omega_1
\end{array}
\right|, \label{two.det} \end{equation} 
where the second determinant is obtained by expanding the first by its last column and
transposing. 
\subsection{} 
On the other hand,
\[
\Delta = \left| 
\binom{j+1}{k} \, \omega^{(j+1-k)} \right|_{0\le k,j\le r}
= \left| \frac{(j+1)!}{k!} \, \bbone_{\{k\le j+1\}} \, \omega_{j+1-k}
\right|_{0\le k,j\le r},   \]
where the characteristic function $\bbone_{\{k\le j+1\}}$ assumes the
value $1$ if $k \le j+1$, and $0$ otherwise. This is the same as the
rightmost determinant in~(\ref{two.det}), hence 
\begin{equation} g_0 \doteq (a_0^{1+\frac{r}{d}})^{r+1} \, \Delta \doteq 
(a_0^{1+\frac{r}{d}})^{r+1} \, \omega_0^{r+2} \,
\nu_{r+1}. \label{g0.final.formula} \end{equation}
Now recall that
\[ \omega_0=\om(0)=a_0^{-\frac{r}{d}}, \]
and
\[
\nu_{r+1}=\frac{1}{(r+1)!} \, \nu^{(r+1)}(0)
=\frac{1}{(r+1)!} \left[\frac{d^{\, r+1}}{d \, t^{\, r+1}} \, f(t)^{\frac{r}{d}} \right]_{t=0}. \]
The exponent of $a_0$ reduces to 
\[
\left(1+\frac{r}{d}\right)(r+1)-\frac{r}{d}(r+2)=r+1-\frac{r}{d}.  \] 
Hence, by comparing~(\ref{h0.formula}) with~(\ref{g0.final.formula}),
we get 
\[ g_0 \doteq h_0,  \] 
which completes the proof of Theorem~\ref{Theorem.HGeq}. \qed 

\bigskip 

\noindent If we keep track of the unwritten scalars in the intermediate stages,
then the connecting relation is 
\begin{equation} g_0 = \left\{\frac{\prod\limits_{i=0}^r \, i! \,
      (d+i-2)!}{\left[ \, r \times (d-2)! \, \right]^{r+1}}
  \, \right\} h_0. \label{formula.krd} \end{equation}
This of course implies a parallel relation between $\Gott_{r,d}$ and
$\Hilb_{r,d}$. 

\subsection{} One can give a formula for the Hilbert covariant directly, 
without constructing its source first. It is merely the homogenised version of the 
formula (\ref{formula.hk}) combined with (\ref{formula.Ediff}), so we 
omit the proof.  Introduce binary variables 
$\uy = \{y_1, y_2\}$. 
\begin{Proposition} \sl 
We have an identity 
\[ \Hilb_{r,d}(\F) \doteq \frac{\F(x_1,x_2)^{r+1-\frac{r}{d}}}{(x_1 \,  y_2 - x_2 \, y_1)^{r+1}} \, 
\left[ \, \left( y_1 \, \frac{\partial}{\partial x_1} + y_2 \,
\frac{\partial}{\partial x_2} \right)^{r+1} \F(x_1,x_2)^{\frac{r}{d}} \right]. \] 
\label{proposition.Hpolar} \end{Proposition} 

\section{The three ideals} 
Let $X = X_{e,d}$ be as in \S\ref{section.defn.rde}, with 
$I_X \subseteq R$ its homogeneous defining ideal. Let $J$ (respectively $\fg$) denote the ideal in $R$
generated by the coefficients of $\Gott_{r,d}$ (respectively all
possible $\Gott_\Psi$). In other words, $\fg$ is the ideal generated
by the maximal minors of a matrix representing the morphism 
$\alpha_\F: S_r \lra S_{r+d-2}$ from \S\ref{definition.alphaF}. There are inclusions 
\[ J \subseteq \fg \subseteq I_X. \] 
The zero locus of each of these ideals is $X$, but depending on the values
of $r$ and $d$, either of these inclusions may be proper. 
Since $I_X$ has nonzero elements in degree $e+1$ (arising from the coefficients of
$\Gott_{e,d}$), we must have a proper containment $\fg \subsetneq
I_X$, if $r$ does not divide $d$. 
\subsection{} \label{section.rnc.ideal} 
Suppose $r=1$, so that $X$ is the rational normal $d$-ic
curve. We have a decomposition 
\[ R_2 \simeq S_2(S_d) \simeq \bigoplus\limits_{n=0}^{\lfloor\frac{d}{2}\rfloor}
\, S_{2d-4n},  \] 
where the summand $S_{2d-4n}$ is spanned by the coefficients of $(\F,\F)_{2n}$. It is classically known
that $I_X$ is minimally generated in degree $2$, and $(I_X)_2 \simeq \bigoplus\limits_{n \ge 1} \, S_{2d-4n}
\subseteq R_2$ (see~\cite{Ch2}). By Proposition \ref{prop.quadG}, we have $\fg =
I_X$. Moreover, $J$ and $\fg$ coincide for $d \le 3$ and differ
afterwards. 
\subsection{} \label{section.g36} 
Assume $r=3, d=6$. One can explicitly calculate the ideal of $X = X_{3,6}$ 
using the following elimination-theoretic technique. Let 
$Q = (q_0,q_1,q_2,q_3 \cb x_1,x_2)^3$, where the $q_i$ are independent
indeterminates. Write 
\[ (a_0, \dots, a_6 \cb x_1,x_2)^6 = Q^2 \] 
and equate the corresponding coefficients on both sides. This
gives expressions $a_i = f_i (q_0, \dots, q_3)$, defining a ring homomorphism 
\[ {\mathfrak f}: R \lra \complex[q_0, \dots, q_3], \qquad
a_i \lra f_i(q_0, \dots, q_3). \] 
Then $I_X$ is the kernel of $\mathfrak f$. We carried out this computation in the
computer algebra system {\sc Macaulay-2} (henceforth {\sc M2}); it
shows that $I_X$ is minimally generated by a $45$-dimensional subspace
of $R_4$. 

In order to determine the piece $(\fg)_4$, we need to list the 
degree $4$ covariants of $\B$. By the Cayley-Sylvester formula, 
\[ S_4(S_4) = S_{16} \oplus S_{12} \oplus S_{10} 
\oplus (S_8 \otimes \complex^2) 
\oplus (S_4 \otimes \complex^2) \oplus  S_0. \] 
It is classically known  (see~\cite[\S 89]{GY}) that each covariant of a generic binary quartic $\B$ is a polynomial in these
fundamental covariants: 
\[ C_{1,4} = \B, \quad C_{2,4} = (\B,\B)_2, \quad C_{2,0} = (\B,\B)_4, \quad
C_{3,6} = (\B, (\B,\B)_2)_1, \quad C_{3,0} = (\B,(\B,\B)_2)_4,  \] 
where $C_{m,q}$ is of degree-order $(m,q)$. Hence, the space of degree $4$ covariants of $\B$ is spanned by 
\[ \begin{array}{lllll} 
\Psi_{4,16} = C_{1,4}^4, & \Psi_{4,12} = C_{1,4}^2 \, C_{2,4}, & \Psi_{4,10} = C_{1,4}
\, C_{3,6}, & \Psi_{4,8}^{(1)} = C_{2,4}^2, & \Psi_{4,8}^{(2)} = C_{1,4}^2 \, C_{2,0}, \\
\Psi_{4,4}^{(1)} = C_{1,4} \, C_{3,0}, & \Psi_{4,4}^{(2)} = C_{2,4} \,
C_{2,0}, & \Psi_{4,0} = C_{2,0}^2. 
\end{array} \] 
We have calculated $\Gott_\Psi$ in each case using the recipe of
\S\ref{GottPsi.formula}. It turns out that the ones coming from $\Psi_{4,16}, \Psi_{4,12},
\Psi_{4,0}$ are nonzero, whereas $\Gott_{\Psi_{4,10}}$ vanishes
identically. Moreover, we have identities 
\[ 6 \, \Gott_{\Psi_{4,8}^{(1)}} = \Gott_{\Psi_{4,8}^{(2)}}, \qquad 
29 \, \Gott_{\Psi_{4,4}^{(1)}} = 36 \, \Gott_{\Psi_{4,4}^{(2)}}, \]
and thus both $\Psi_{(4,8)}^{(i)}$ lead to the same G{\"o}ttingen covariant (up to a
scalar), and similarly for $\Psi_{(4,4)}^{(i)}$. Hence 
\[ (\fg)_4 \simeq S_{16} \oplus S_{12} \oplus S_8 \oplus S_4 \oplus
S_0, \] 
which is exactly $45$-dimensional;  this forces $\fg = I_X$. 

\subsection{} \label{section.g26} 
Assume $r=2$, and $d$ even. Now \cite[Theorem 7.2]{AC1} says that
$I_X$ is minimally generated by cubic forms, and its generators are explicitly
described there. If $d=4$, then $(I_X)_3 \simeq S_6$, with the only
piece coming from $\Gott_{2,4}$. If $d=6$, then 
\[ (I_X)_3 \simeq S_{12} \oplus S_8 \oplus S_6.  \] 
The three summands are respectively generated by the coefficients
of: 
\[ \Phi_{3,12} = (\F^2,\F)_3, \quad \Phi_{3,8} = (\F^2,\F)_5, \quad
\Phi_{3,6}  = 33 \, (\F^2,\F)_6 - 250 \, (\F,(\F,\F)_2)_4.  \] 
Now, following the recipe of \S \ref{GottPsi.formula}, one finds that 
\[ \Gott_{\B^3} \doteq \Phi_{3,12}, \quad 
\Gott_{\B \, (\B,\B)_2} \doteq \Phi_{3,8}, \quad 
\Gott_{(\B, (\B,\B)_2)_1} \doteq \Phi_{3,6};   \] 
and hence $\fg = I_X$ once again. 

\medskip 

We have calculated several such examples, which suggest the
following pair of conjectures: 
\begin{Conjecture} \sl 
Assume that $r$ divides $d$. Then 
\begin{itemize} 
\item[(c1)] 
the ideal $I_X$ is always minimally generated in degree $r+1$, and 
\item[(c2)] $\fg = I_X$. 
\end{itemize} 
\end{Conjecture} 

At least for $r=1, 2$, something much stronger than (c1) is true;
namely $I_X$ has Castelnuovo regularity $r+1$, and its graded minimal
resolution is linear (see~\cite[Theorem 1.4]{AC1}). We do not know of
a counterexample to this when $r> 2$. 

Referring to the diagram at the beginning of \S\ref{section.Gott.cov}, note that the ideal of the Veronese
embedding is generated by quadrics; but the projection $\pi$ (which
implicitly involves elimination theory) will tend to increase the degrees of the
defining equations of its image. 

\subsection{The saturation of $J$} \label{proof.Jsaturation} 
In this section we will prove Theorem~\ref{Theorem.Jsaturation}. We
want to show that the ideal $J$ defines $X$ scheme-theoretically when $r$ divides $d$, that is to say, 
\[ \Proj R/I_X \lra \Proj R/J \] 
is an isomorphism of schemes.  The following example
should convey the essential idea behind the proof. 

Assume $r=2, d=6$. Write $t_i = a_i/a_0$ for $1 \le i \le 6$. Let 
$A = \complex[t_1, \dots, t_6]$, and consider the corresponding degree
zero localisation $\fra = (J_{a_0})_0 \subseteq A$. The zero
locus of $\fra$ is $X \setminus \{a_0=0\} \simeq \Aff^2$. Since the
question is local on $X$, it would suffice to show that $A/\fra$ is isomorphic to a polynomial algebra 
$\complex[v_1,v_2]$. 

Now $\Hilb_{2,d} \doteq (\F, (\F, \F)_2)_1$, and we have explicitly
written down its first few terms in \S\ref{ex.isobaric}. Note that the
monomial $a_0^2 \, a_3$ occurs in its source, and similarly $a_0^2 \, a_4, a_0^2 \, a_5, a_0^2
\, a_6$ occur in the successive coefficients. Hence, modulo $\fra$, we have identities of
the form 
\[ t_k = \text{a polynomial expression in $t_1, t_2, \dots, t_{k-1}$,}
\quad \text{for $3\le k \le 6$}. \] 
Thus we have a surjective ring morphism 
\[ \complex[v_1,v_2] \lra A/\fra, \quad v_i \lra t_i.  \] 
Since $\text{Krull-dim} \, A/\fra=2$, this must be an isomorphism. 

\medskip 

For the general case, write $\Hilb_{r,d}$ as in~(\ref{Hilb.cov.definition}),
and recall that $h_{k-(r+1)}$ is isobaric of weight $k$. 
\begin{Lemma} \sl The coefficient of $a_0^r \, a_k$ in $h_{k-(r+1)}$ is
  nonzero for $r+1 \le k \le d$. 
\end{Lemma} 
\demo 
The monomial $a_0^r \, a_{r+1}$ can appear in $h_0$ only by one route,
namely by applying the sequence 
\[ \left[ (d-r) \, a_{r+1} \, \frac{\partial}{\partial a_r} \right]
\circ \dots \left[(d-1) \, a_2 \, \frac{\partial}{\partial a_1}
\right] \circ \left[ d \, a_1 \,\frac{\partial}{\partial a_0} \right] \] 
to $a_0^{\frac{r}{d}}$, and then multiplying by
$a_0^{r+1-\frac{r}{d}}$. Hence its coefficient is nonzero. 
Now $a_0^r \, a_k$ can appear in $h_{k-(r+1)} \doteq E_+ \, h_{k-1-(r+1)}$ only by 
applying $(d-k+1) \, a_k \, \frac{\partial}{\partial a_{k-1}}$ to
$a_0^r \, a_{k-1}$, so we are done by induction. \qed 

\medskip 

We can always change co-ordinates such that $a_0 \neq 0$ at any
given point of $X$. Write $t_i = a_i/a_0$ and $\fra < A = \complex[t_1,
\dots, t_d]$ as above. By the lemma, each of $t_{r+1}, \dots,
t_d$ is a polynomial in $t_1, \dots, t_r$ modulo $\fra$. This gives a bijection 
\[ \complex[v_1,\dots, v_r] \lra A/\fra, \quad v_i \lra t_i,   \] 
which shows that the scheme $\Proj R/J$ is locally isomorphic
to the affine space $\Aff^r$, and hence $J_\text{sat}  = I_X$. 
This completes the proof of Theorem~\ref{Theorem.Jsaturation}. \qed 

\subsection{} It follows that $J$ and $I_X$ coincide in sufficiently
large degrees. Let $\SI(r,d)$ denote the saturation index of $J$, 
namely it is the smallest integer $m_0$ such that
\[ J_m = (I_X)_m \qquad \text{for all \, $m \ge m_0$}. \] 
It would be of interest to have a bound on this quantity in either
direction. It is proved in~\cite{Ch2} that 
\[ \frac{1}{d-2} \sqrt{\frac{(d-1)(d^2-2)}{2}} \le \SI(1,d) \le
d+2;  \] 
but those techniques do not seem to generalise readily to the case
$r>1$. We have obtained the following few values by explicit
calculations in {\sc M2:}  
\[ \begin{array}{lllll} 
\SI(2,4) = 3, & \SI(2,6) = 7, & \SI(2,8) = 9, & \SI(2,10) = 9, & \SI(2,12) = 10, \\ 
\SI(3,6) = 9, & \SI(3,9) = 11, & {\mathfrak S}(4,8) = 13. 
\end{array} \] 
A similar (but larger) table for $r=1$ is given in \cite{Ch2}, where the value of
$\SI$ is related to transvectant identities involving the Hessian. 

\subsection{} Suppose $e_i = \gcd (r_i,d)$ for $i=1,2$. Then $X_{e_1,d}
\subseteq X_{e_2,d}$ exactly when $e_1 \, | \, e_2$. 
However, the containment relations between the ideals $J_{r_i,d}$ are
not altogether obvious. For $J_{r_1,d} \supseteq J_{r_2,d}$ 
to be true, it is necessary that $r_1 \le r_2$ and $e_1 \, |
\, e_2$, but these conditions are not sufficient. For
instance, we have obtained the following miscellaneous data by calculating these ideals in {\sc M2}: 
\[ \begin{array}{llll} 
J_{2,5} \not\supseteq J_{3,5}, & J_{3,5} \not\supseteq J_{4,5}, &
J_{2,5} \supseteq J_{4,5}, & J_{4,5} \not\supseteq J_{6,5}, \\ 
J_{2,4} \supseteq J_{6,4}, & J_{6,4} \supseteq
J_{10,4};   \end{array} \] 
which at least shows that the general pattern is not so easily
guessed. Nevertheless, we have the following modest result: 
\begin{Proposition} \sl 
There are inclusions $J_{1,d} \supseteq J_{r,d}$ for arbitrary $d$, and $r=2,3,4$. 
\label{prop.containmentJ} \end{Proposition} 
\demo It is clear from the formula for a transvectant (see
\S\ref{section.trans}), that if the coefficients of $A$
belong to an ideal, then all the coefficients of $(A,B)_k$
also belong to this ideal. Hence, given any covariants $\Phi_1, \dots,
\Phi_n$ of $\F$, all the coefficients of any transvectant of the form 
\[ (\dots ((\Gott_{1,d}, \Phi_1)_{k_1}, \Phi_2)_{k_2}, \dots , \Phi_n)_{k_n} \] 
are in $J_{1,d}$. Thus the result would follow if we could obtain $\Gott_{r,d}$ as a linear
combination of such expressions. 

Observe the formulae in (\ref{Gott.lowr}). It is clear that
$\Gott_{2,d}$ is itself such an expression. Let $r=4$, then this is also true
of the first term in $\Gott_{4,d}$. Now the so-called Gordan syzygies give relations between cubic covariants
of $\F$. In particular, the syzygy which is written as
$\left( \begin{array}{ccc} \F & \F & \F \\ d & d & d \\ 0 & 1 &
    4 \end{array} \right)$ in the notation of \cite[Ch.~IV]{GY}, gives an
identity 
\[ (\F, (\F, \F)_4)_1 = \frac{2 \, (2d-5)}{d-4} \, (\F, (\F, \F)_2)_3,  \] 
for any $d \ge 5$. Hence the same follows for the second
term. (If $d \le 4$, then the second term is identically zero.) This proves the result for $r=4$. 

The argument for $r=3$ is similar. The first term in $\Gott_{3,d}$
is already of the required form. Moreover, we have an identity 
\[ \F^2 \, (\F,\F)_4 = \frac{d \, (2d-5)}{(d-3) \, (2d-1)} \, (\F,\F)_2^2
+ \frac{2 \, (2d-5)}{d-3} \, (\F^2,(\F,\F)_2)_2,  \] 
for $d \ge 4$. (This can be shown by a routine but tedious symbolic calculation as in
\cite[Ch.~IV-V]{GY}.) Hence the same is true of the second term, which completes the proof. \qed 

\smallskip 

Since the argument depends on specific features of these formulae, 
it seems unlikely that this technique will generalise
substantially. Even so, we suspect that the proposition may well be true of all $r$. 

\subsection{The twisted cubic curve} \label{section.twistedcubic} 
Assume $d=3$, and $r$ arbitrary (but
not divisible by $3$). Then $X \subseteq \P^3$ is the twisted cubic curve. Since $\B$ is a linear form, the only
possibility for $\Psi$ is $\B^{r+1}$, hence $J = \fg$. It follows that 
the Hilbert-Burch complex (see~\cite[\S 20]{Eisenbud1}) of $\alpha_\F$
gives a resolution
\[ 0 \la R/J \la R \la R(-r-1) \otimes S_{r+1} \la R(-r-2) \otimes S_r
\la 0. \] 
Its first syzygy shows that we have an identity $(\Gott_{r,3},
\F)_2=0$. (The correspondence between syzygies and transvectant
identities is discussed in \cite[\S 4]{Ch1}.) The scheme $\Proj R/J$
has degree $\binom{r+2}{2}$, that is to say, it is a nonreduced
$\frac{(r+1) \, (r+2)}{6}$--fold structure on $X$ for $r>1$. 

We have $\sqrt{J_{r,3}} = I_X$ for any $r$. Some experimental calculations in {\sc M2} suggest the
following narrow but interesting conjecture: 
\begin{Conjecture} \sl There is an inclusion 
$(I_X)^r \subseteq J_{r,3}$, and moreover $r$ is the smallest such power. 
\end{Conjecture}

This problem is related to identities between the covariants of a
generic cubic form. For instance, we have an identity 
\[ \Gott_{1,3}^2 = - \frac{1}{2} \, (\F,\Gott_{2,3})_1,  \] 
which can be verified by a direct symbolic computation. This
immediately shows that $(I_X)^2 \subseteq J_{2,3}$. (Compare the
argument of Proposition~\ref{prop.containmentJ} above.) 

In general, if $r$ and $d$ are coprime, then $\fg$ is a perfect ideal of height $d-1$,
which is resolved by the Eagon-Northcott complex (see \cite[Appendix 2]{Eisenbud1}) of $\alpha_\F$. 
By the Porteous formula (see~\cite[Ch.~II, \S 4]{ACGH}), the scheme $\Proj R/\fg$ supported on the
rational normal $d$-ic curve has degree $\binom{r+d-1}{d-1}$. 
\section{The Clebsch transfer principle} \label{section.Clebsch}
In this section we generalise the G{\"o}ttingen covariants to $n$-ary
forms. 
\subsection{} 
Let $W$ be an $n$-dimensional complex vector space with basis $\ux = \{x_1,
\dots, x_n\}$, and a natural action of the group $SL(W)$. 
Given an $n$-tuple of nonnegative integers $I = (i_1, \dots, i_n)$ adding up to
$d$, let 
\[ \binom{d}{I} = \frac{d!}{\prod\limits_k \, i_k!}, \qquad x^I =
\prod\limits x_k^{i_k}.   \] 
We write a generic form of order $d$ in the $\ux$ as 
\[ \Gamma = \sum\limits_I \binom{d}{I} \, a_I \, x^I,  \] 
where the $a_I$ are independent indeterminates. As in the binary case,
the $\{a_I\}$ can be seen as forming a basis of $S_d \, W^*$. Define the symmetric algebra 
\[ \cA = \bigoplus\limits_{m \ge 0} \, S_m(S_d \, W^*) = \complex[\{a_I\}] \] 
so that $\Proj \cA = \P S_d \simeq \P^{\binom{d+n-1}{d}-1}$ is the space of
$n$-ary $d$-ics. 

\subsection{} 
Each irreducible representation of $SL(W)$ is a
Schur module of the form $S_\lambda = S_\lambda \, W$, where $\lambda$ is a
partition with at most $n-1$ parts (see~\cite[\S 15]{FH}). Moreover,
we have an isomorphism 
\[ S_{(\lambda_1, \lambda_2, \dots, \lambda_{n-1},0)}  \, W \simeq 
S_{(\lambda_1, \lambda_1-\lambda_{n-1}, \dots, \lambda_1 -
  \lambda_2,0)}  \, W^*. \] 
An inclusion $S_\lambda \, W^* \subseteq \cA_m$ corresponds to a morphism 
\[ S_0 \hookrightarrow \cA_m \otimes S_\lambda \, W,  \] 
then the image of $1 \in S_0$ will be called a concomitant of $\Gamma$ of degree
$m$, and type $\lambda$. 

\subsection{} \label{section.clebsch33} 
In the case of ternary forms, for $\lambda = (\lambda_1, \lambda_2)$, 
we have an embedding (see~\cite[\S 15]{FH}) 
\[ S_\lambda \subseteq S_{\lambda_2}(\wedge^2 \, W) \otimes S_{\lambda_1
  - \lambda_2}. \] 
Using the basis $u_1 = x_1 \wedge x_2, \, u_2 = x_2 \wedge x_3,
\, u_3 = x_3 \wedge x_1$ for $\wedge^2 \, W \simeq W^*$, we can write the concomitant as a form of degree
$m$ in the $a_I$, degree $\lambda_1 - \lambda_2$ in $\ux$, and
degree $\lambda_2$ in $\uu$. 
For instance, assume $m=2,d=3$. We have a plethysm decomposition 
$S_2(S_3^*) \simeq S_6^* \oplus S_{4,2}^*$, and hence (up to a scalar) a unique morphism 
\[ S_0 \hookrightarrow \cA_2 \otimes S_{4,2}. \] 
If we symbolically write $\Gamma = a_\ux^3 = b_\ux^3$, then this
concomitant is $(a \, b \, u)^2 \, a_\ux \, b_\ux$. We refer the reader
to~\cite[Ch.~XII]{GY} or~\cite{Littlewood} for 
the symbolic calculus of $n$-ary forms and their concomitants. 

\subsection{} \label{section.clebsch34} 
The `Clebsch transfer principle' is a type of construction used to lift a binary covariant to a concomitant of $n$-ary
forms in a geometrically natural way. As such, it comes in many
flavours depending on the specifics of the geometric situation in
play. (See~\cite[\S 4]{Briand}, \cite[\S 3.4.2]{Dolgachev2} or \cite[\S 215]{GY} for variant descriptions
of this principle.) Clebsch's own statement of this technique may be found
in~\cite[p.~28]{Clebsch}, but Cayley and Salmon seem to have been
aware of it earlier (see~\cite[p.~28]{SalmonCD}). 

The following example should convey an idea of how the transfer principle is used. 
Let $n=3$ and $d=4$, so that $\P S_4 \simeq \P^{14}$ is the space of
quartic plane curves. Let $Z \subset \P^{14}$ be the $5$-dimensional
subvariety of double conics, i.e., 
\[ Z = \left\{  [\Gamma] \in \P S_4: \Gamma = Q^2 \; \, \text{for some
  ternary quadratic $Q$} \right\}. \] 
A line $L$ in the plane $\P W^* \simeq \P^2$ will intersect a general 
quartic curve $\Gamma(x_1,x_2,x_3)=0$ in four points, which become 
two double points when $\Gamma \in Z$. With the identification $L
\simeq \P^1$, let $\Gamma|_L$ denote the `restriction' of $\Gamma$ to
$L$, regarded as a binary quartic form. Hence the `function' 
\[ L \lra \Gott_{2,4}(\Gamma|_L) \] 
should vanish identically when $\Gamma \in Z$. 

In order to make this precise, write $p = [p_1, p_2, p_3], q =
[q_1,q_2,q_3]$, where $p_i,q_i$ are indeterminates. We think of a generic $L$ as spanned by the points $p,q \in
\P^2$, and thus $L$ has line co-ordinates 
\[ u_1 = \left| \begin{array}{cc} p_1 & q_1 \\ p_2 & q_2 \end{array}
\right|, \quad 
u_2 = \left| \begin{array}{cc} p_2 & q_2 \\ p_3 & q_3 \end{array} \right|,
\quad 
u_3 = \left| \begin{array}{cc} p_3 & q_3 \\ p_1 & q_1 \end{array}
\right|. \] 
Introduce binary variables $\lambda = \{\lambda_1, \lambda_2 \}$,
and substitute $x_i = \lambda_1 \, p_i + \lambda_2 \, q_i$ in
$\Gamma$ to get a form $\Theta$ (which represents the restriction). 
Now evaluate $\Gott_{2,4}$ on $\Theta$ by regarding the
latter as a binary form in the $\lambda$; then the final result is the
required lift $\tGott_{2,4}$. The actual symbolic calculation proceeds as follows. 
Let 
\[ \Gamma = a_\ux^4 = b_\ux^4 = c_\ux^4, \] 
where $a_\ux = a_1 \, x_1 + a_2 \, x_2 + a_3 \, x_3$ etc. After
substitution, $a_\ux$ becomes $\lambda_1 \, a_p + \lambda_2 \, a_q$, which we rewrite as
\[ \alpha_\lambda = \alpha_1 \, \lambda_1 + \alpha_2 \, \lambda_2 \quad 
\text{where $\alpha_1 = a_p, \, \alpha_2 = a_q$,} \] 
and similarly $b_\ux = \beta_\lambda, c_\ux = \gamma_\lambda$. 
Thus $\Theta = \alpha_\lambda^4 = \beta_\lambda^4 = \gamma_\lambda^4$. 
Recall from \S\ref{Gott.quadratic} that 
\begin{equation} 
\Gott_{2,d}(\Theta) \doteq (\Theta, (\Theta,\Theta)_2)_1 \doteq 
(\alpha \, \beta)^2 \, (\alpha \, \gamma) \, \alpha_\lambda \, \beta_\lambda^2
\, \gamma_\lambda^3. \label{Gott.24.Theta} \end{equation} 
Now 
\[ (\alpha \, \beta) = \left| \begin{array}{cc} \alpha_1 & \alpha_2 \\
  \beta_1 & \beta_2 \end{array} \right| = a_p \, b_q - b_p \, a_q = 
\left| \begin{array}{ccc} a_1 & a_2 & a_3 \\ b_1 & b_2 & b_3 \\ u_1 &
    u_2 & u_3 \end{array} \right| = (a \, b \, u),  \] 
and similarly for the bracket factor $(\alpha \, \gamma)$. Hence we
arrive at the expression 
\begin{equation} \tGott_{2,4}(\Gamma) = (a \, b \, u)^2 \, (a \, c \, u) \, a_\ux \,
b_\ux^2 \, c_\ux^3,  \label{tGott.24} \end{equation}
which is a concomitant of degree $3$ and type $(9,3)$. We have the property
\[ [\Gamma] \in Z \iff \tGott_{2,4}(\Gamma) \; \text{vanishes identically as a
polynomial in $\ux, \uu$}. \] 
The implication $\Rightarrow$ follows by construction. The converse says 
that if $\Gamma=0$ were not a double conic, then a line could be found which does not intersect it in two
double points. This is clear on geometric grounds. 

\subsection{} The case of a general G{\"o}ttingen covariant is
similar. Assume that $\Gott_\Psi$ is of degree $r+1$, order $q$, and
weight $w=\frac{(r+1) \, d-q}{2}$. Let 
$p=[  p_1, \dots, p_n], \, q=[ q_1, \dots, q_n]$, substitute 
\begin{equation} x_i = \lambda_1 \, p_i + \lambda_2 \, q_i, \quad  (1 \le i \le n), 
\label{subs.xi} \end{equation} 
into $\Gamma$, and evaluate $\Gott_\Psi$ on the new binary form in the $\lambda$
variables. The resulting concomitant $\tGott_\Psi$ is of degree $r+1$
and type $(q+w,w)$. If $\Gamma = G^\mu$, then $\tGott_\Psi(\Gamma)$ vanishes identically for the same
reason as above. 

If $\Gott_\Psi$ is written as a symbolic expression in $r+1$ binary letters
$a,b,\dots$ and their brackets $(a \, b)$ etc., then $\tGott_\Psi$ is
obtained by simply treating them as $n$-ary letters and replacing the corresponding brackets by $(a \, b \, u)$
etc. This follows immediately by tracing the passage from~(\ref{Gott.24.Theta}) to~(\ref{tGott.24}). 
In particular, the concomitant in \S\ref{section.clebsch33} is the
Clebsch transfer of the Hessian of a binary cubic. The formal symbolic
expression for $\tGott_\Psi$ does not depend on $n$, although of
course, its interpretation does. 

\begin{Theorem} \sl 
Let $\Gamma$ be an $n$-ary $d$-ic. Then $\tGott_{r,d}(\Gamma)$ is
identically zero, if and only if $\Gamma = G^\mu$ for some $n$-ic $G$ of
order $e$. 
\label{theorem.clebsch} \end{Theorem} 
\demo The `if' part follows from the discussion above. Let $\Gamma =
\prod\limits H_i^{\nu_i}$ be the prime decomposition, where $H_i$ is
an irreducible form of degree $c_i$. If $\Gamma$ cannot be written as $G^\mu$, then at
least one $\nu_i$ is not divisible by $\mu$. A general line $L$ will
intersect each hypersurface $H_i=0$ in $c_i$ distinct
points. Altogether $L$ intersects $\Gamma=0$ in $c_1+ c_2 + \dots$
points, at least one of which occurs with multiplicity not divisible by
$\mu$. Thus $\tGott_{r,d}(\Gamma)$ will not vanish if the $\uu$ variables in
it are specialised to the Pl{\"u}cker co-ordinates of a general $L$. \qed 

\subsection{} This is a continuation of \S\ref{section.clebsch34}. We have calculated the ideal of $Z$ 
using a procedure similar to the one in \S\ref{section.g36}, and it turns
out that $I_Z$ is minimally generated by a $218$-dimensional space of forms in
degree $3$. We have a plethysm decomposition 
\[ \cA_3 = S_3(S_4^*) \simeq S_{12}^* \oplus S_{10,2}^* \oplus S_{9,3}^*
\oplus S_{8,4}^* \oplus S_6^* \oplus S_{6,3}^* \oplus S_{6,6}^* \oplus
S_{4,2}^* \oplus S_0^*,  \] 
where the summands are of respective dimensions 
\[ 91, \, 162, \, 154, \, 125, \, 28, \, 64, \, 28, \, 27, \, 1. \] 
Now $(I_Z)_3$ is a subdirect sum of the above, and we
already know that $S_{9,3}^*$ is one of its pieces. This forces 
$(I_Z)_3 \simeq S_{9,3}^* \oplus S_{6,3}^*$ on dimensional
grounds. Hence there is a concomitant of type $(6,3)$ vanishing on
$Z$. We have checked by a direct calculation that it can be
written as 
\[ (a \, b \, c) \, (a \, b \, u)^2 \, (a \, c \, u) \, b_\ux \,
c_\ux^2. \] 
In fact, all that needs to be checked is that this symbolic expression
is not identically zero, which can be done by specialising
$\Gamma$. This suffices, since we have up to scalar only one concomitant of this type in degree $3$. 

Recall from \S\ref{section.g26} that for $r=2, d=4$, that there are 
no G{\"o}ttingen covariants other than $\Gott_{2,4}$. Hence we have
found a concomitant vanishing on $Z$ which is not the Clebsch transfer
of any binary covariant. 

Let $J \subseteq \cA$ denote the ideal generated by the coefficients of
$\tGott_{2,4}$. We have checked using {\sc M2} that
the saturation of $J$ is $I_Z$, and moreover the two ideals coincide
in degrees $\ge 7$. But in general, we do not know whether there is an
analogue of Theorem \ref{Theorem.Jsaturation} in the $n$-ary case. 

\subsection{} We end with an example which is at least a pleasing
curiosity. Assume that $\Gamma=0$ is a \emph{nonsingular} plane quartic
curve. A line $L \subset \P^2$ with co-ordinates $[u_1,u_2,u_3]$
passes through the points $p = [u_3,0,-u_1], q = [u_2,-u_1,0]$, 
and moreover these points are distinct (and well-defined) when $u_1
\neq 0$. Now make substitutions into $\tGott_{2,4}(\Gamma)$ as in~(\ref{subs.xi}) to get a binary sextic
$\cE_1(\lambda)$; it represents the binary form $\Gott_{2,4}$ as
living on $L \simeq \P^1$. (This is no longer correct if $u_1=0$, hence
in order to avoid spurious solutions we also 
need to consider the forms $\cE_2(\lambda), \cE_3(\lambda)$ similarly obtained from 
\[ p = [0,u_3,-u_2], \; q=[u_2,-u_1,0], \qquad p = [0,u_3,-u_2], \;
q=[u_3,0,-u_1].)  \] 
Now all the $\cE_i(\lambda)$ are identically
zero exactly when $\{\Gamma=0\} \cap L$ represents two double points, i.e., when $L$ is a
bitangent to the curve defined by $\Gamma$. Let $B = \complex[u_1,u_2,u_3]$ denote the
co-ordinate ring of the dual plane, and ${\mathfrak b}_\Gamma
\subseteq B$ the ideal generated by the coefficients of all the monomials
in $\lambda$ for $\cE_i(\lambda), i=1,2,3$. Then the zero locus of ${\mathfrak
  b}_\Gamma$ is the set of $28$ points  (see~\cite[Ch.~6]{Dolgachev2})
corresponding to the bitangents of the curve. 
We have verified in {\sc M2} that ${\mathfrak b}_\Gamma$ is not
saturated, but its saturation has resolution 
\[ 0 \la B/({\mathfrak b}_\Gamma)_\text{sat} \la B \la B(-7)^8 \la B(-8)^7 \la 0,  \] 
which is characteristic of $28$ general points in the plane
(see~\cite[Ch.~3]{Eisenbud2}). In much the same way, the concomitant
in \S\ref{section.clebsch33} can be used to give equations for the
$9$ inflexional tangents of a nonsingular plane cubic curve. 

\bigskip 

{\small {\sc Acknowledgements:} The second author was funded by NSERC, Canada. We are thankful to Daniel Grayson and
Michael Stillman (the authors of {\sc Macaulay-2}). We have used John Stembridge's `SF'
package for {\sc Maple} for calculating plethysm decompositions.}

\medskip 

\centerline{---} 

\vspace{1cm} 

\parbox{7cm} 
{Abdelmalek Abdesselam \\
Department of Mathematics, \\
University of Virginia, \\
P. O. Box 400137, \\
Charlottesville, VA 22904-4137, \\
USA.\\
{\tt malek@virginia.edu}} 
\hfill 
\parbox{7cm} 
{Jaydeep Chipalkatti \\ 
Department of Mathematics, \\ 
Machray Hall, \\ 
University of Manitoba, \\ 
Winnipeg, MB R3T 2N2, \\ Canada. \\ 
{\tt chipalka@cc.umanitoba.ca}}

\end{document}